\documentclass{article}%
\usepackage{amsfonts}
\usepackage{amssymb}
\usepackage{amsmath}
\usepackage{graphicx}%
\providecommand{\U}[1]{\protect\rule{.1in}{.1in}}
\pdfoutput=1

\newtheorem{corollary}{Corollary}

\newtheorem{proposition}{Proposition}

\begin{document}

\title{The number of direct-sum decompositions\\of a finite vector space}
\author{David Ellerman\\University of California at Riverside}
\maketitle

\begin{abstract}
\noindent The theory of $q$-analogs develops many combinatorial formulas for
finite vector spaces over a finite field with $q$ elements--all in analogy
with formulas for finite sets (which are the special case of $q=1$). A
direct-sum decomposition of a finite vector space is the vector space analogue
of a set partition. This paper uses elementary methods to develop the formulas
for the number of direct-sum decompositions that are the $q$-analogs of the
formulas for: (1) the number of set partitions with a given number partition
signature; (2) the number of set partitions of an $n$-element set with $m$
blocks (the Stirling numbers of the second kind); and (3) for the total number
of set partitions of an $n$-element set (the Bell numbers).

\end{abstract}
\tableofcontents

\section{Reviewing q-analogs: From sets to vector spaces}

The theory of $q$-analogs shows how many "classical" combinatorial formulas
for finite sets can be extended to finite vector spaces where $q$ is the
cardinality of the finite base field, i.e., $q=p^{n}$, a power of a prime.

The natural number $n$ is replaced by:

\begin{center}
$\left[  n\right]  _{q}=\frac{q^{n}-1}{q-1}=1+q+q^{2}+...+q^{n-1}$
\end{center}

\noindent so as $q\rightarrow1$, then $\left[  n\right]  _{q}\rightarrow n$ in
the passage from vector spaces to sets. The factorial $n!$ is replaced, in the
$q$-analog

\begin{center}
$\left[  n\right]  _{q}!=\left[  n\right]  _{q}\left[  n-1\right]
_{q}...\left[  1\right]  _{q}$
\end{center}

\noindent where $\left[  1\right]  _{q}=\left[  0\right]  _{q}=1$.

To obtain the Gaussian binomial coefficients we calculate with ordered bases
of a $k$-dimensional subspace of an $n$-dimensional vector space over the
finite field $GF\left(  q\right)  $ with $q$ elements. There are $q^{n}$
elements in the space so the first choice for a basis vector has $\left(
q^{n}-1\right)  $ (excluding $0$) possibilities, and since that vector
generated a subspace of dimension $q$, the choice of the second basis vector
is limited to $\left(  q^{n}-q\right)  $ elements, and so forth. Thus:

\begin{center}
$\left(  q^{n}-1\right)  \left(  q^{n}-q\right)  \left(  q^{n}-q^{2}\right)
...\left(  q^{n}-q^{k-1}\right)  $

$=\left(  q^{n}-1\right)  q^{1}\left(  q^{n-1}-1\right)  q^{2}\left(
q^{n-1}-1\right)  ...q^{k-1}\left(  q^{n-k+1}-1\right)  $

$=\frac{\left[  n\right]  _{q}!}{\left[  n-k\right]  _{q}!}q^{\left(
1+2+...+\left(  k-1\right)  \right)  }=\frac{\left[  n\right]  _{q}!}{\left[
n-k\right]  _{q}!}q^{k\left(  k-1\right)  /2}$.

Number of ordered bases for a $k$-dimensional subspace in an $n$-dimensional space.
\end{center}

\noindent But for a space of dimension $k$, the number of ordered bases are:

\begin{center}
$\left(  q^{k}-1\right)  \left(  q^{k}-q\right)  \left(  q^{k}-q^{2}\right)
...\left(  q^{k}-q^{k-1}\right)  $

$=\left(  q^{k}-1\right)  q^{1}\left(  q^{k-1}-1\right)  q^{2}\left(
q^{k-1}-1\right)  ...q^{k-1}\left(  q^{k-k+1}-1\right)  $

$=\left[  k\right]  _{q}!q^{k\left(  k-1\right)  /2}$

Number of ordered bases for a $k$-dimensional space.
\end{center}

Thus the number of subspaces of dimension $k$ is the ratio:

\begin{center}
$\binom{n}{k}_{q}=\frac{\left[  n\right]  _{q}!q^{k\left(  k-1\right)  /2}%
}{\left[  n-k\right]  _{q}!\left[  k\right]  _{q}!q^{k\left(  k-1\right)  /2}%
}=\frac{\left[  n\right]  _{q}!}{\left[  n-k\right]  _{q}!\left[  k\right]
_{q}!}$

Gaussian binomial coefficient
\end{center}

\noindent where $\binom{n}{k}_{q}\rightarrow\binom{n}{k}$ as $q\rightarrow1$,
i.e., the number of $k$-dimensional subspaces $\rightarrow$ number of
$k$-element subsets. Many classical identities for binomial coefficients
generalize to Gaussian binomial coefficients \cite{goldman-rota:foundations4}.

\section{Counting partitions of finite sets and vector spaces}

\subsection{The direct formulas for counting partitions of finite sets}

Using sophisticated techniques, the direct-sum decompositions of a finite
vector space over $GF\left(  q\right)  $ have been enumerated in the sense of
giving the exponential generating function for the numbers
(\cite{bender-goldman}; \cite{stanley:expstructures}). Our goal is to derive
by elementary methods the formulas to enumerate these and some related
direct-sum decompositions.

Two subspaces of a vector space are said to be \textit{disjoint} if their
intersection is the zero subspace $0$. A \textit{direct-sum decomposition}
(DSD) of a finite-dimensional vector space $V$ over a base field $F$ is a set
of (nonzero) pair-wise disjoint subspaces, called \textit{blocks} (as with
partitions), $\left\{  V_{i}\right\}  _{i=1,...,m}$ that span the space. Then
each vector $v\in V$ has a unique expression $v=\sum_{i=1}^{m}v_{i}$ with each
$v_{i}\in V_{i}$. Since a direct-sum decomposition can be seen as the
vector-space version of a set partition, we begin with counting the number of
partitions on a set.

Each set partition $\left\{  B_{1},...,B_{m}\right\}  $ of an $n$-element set
has a "type" or "signature" number partition giving the cardinality of the
blocks where they might be presented in nondecreasing order which we can
assume to be: $\left(  \left\vert B_{1}\right\vert ,\left\vert B_{2}%
\right\vert ,...,\left\vert B_{m}\right\vert \right)  $ which is a number
partition of $n$. For our purposes, there is another way to present number
partitions, the \textit{part-count representation}, where $a_{k}$ is the
number of times the integer $k$ occurs in the number partition (and $a_{k}=0$
if $k$ does not appear) so that:

\begin{center}
$a_{1}1+a_{2}2+...+a_{n}n=\sum_{k=1}^{n}a_{k}k=n$.

Part-count representation of number partitions keeping track of repetitions.
\end{center}

Each set partition $\left\{  B_{1},...,B_{m}\right\}  $ of an $n$-element set
has a part-count signature $a_{1},...,a_{n}$, and then there is a "classical"
formula for the number of partitions with that signature (\cite[p.
215]{andrews:partitions}; \cite[p. 427]{knuth:vol4a}).

\begin{proposition}
The number of set partitions for the given signature: $a_{1},...,a_{n}$ where
$\sum_{k=1}^{n}a_{k}k=n$ is:
\end{proposition}

\begin{center}
$\frac{n!}{a_{1}!a_{2}!...a_{n}!\left(  1!\right)  ^{a_{1}}\left(  2!\right)
^{a_{2}}...\left(  n!\right)  ^{a_{n}}}$.

.
\end{center}

\noindent Proof: Suppose we count the number of set partitions $\left\{
B_{1},...,B_{m}\right\}  $ of an $n$-element set when the blocks have the
given cardinalities: $n_{j}=\left\vert B_{j}\right\vert $ for $j=1,...,m$ so
$\sum_{j=1}^{m}n_{j}=n$. The first block $B_{1}$ can be chosen in $\binom
{n}{n_{1}}$ ways, the second block in $\binom{n-n_{1}}{n_{2}}$ ways and so
forth, so the total number of ways is:

\begin{center}
$\binom{n}{n_{1}}\binom{n-n_{1}}{n_{2}}...\binom{n-n_{1}-...-n_{m-1}}{n_{m}%
}=\frac{n!}{n_{1}!\left(  n-n_{1}\right)  !}\frac{\left(  n-n_{1}\right)
!}{n_{2}!\left(  n-n_{1}-n_{2}\right)  !}...\frac{\left(  n-n_{1}%
-...-n_{m-1}\right)  !}{n_{m}!\left(  n-n_{1}-...-n_{m}\right)  !}$

$=\frac{n!}{n_{1}!...n_{m}!}=\binom{n}{n_{1},...,n_{m}}$
\end{center}

\noindent the multinomial coefficient. This formula can then be restated in
terms of the part-count signature $a_{1},...,a_{n}$ where $\sum_{k=1}^{n}%
a_{k}k=n$ as: $\frac{n!}{\left(  1!\right)  ^{a_{1}}\left(  2!\right)
^{a_{2}}...\left(  n!\right)  ^{a_{n}}}$. But that overcounts since the
$a_{k}$ blocks of size $k$ can be permuted without changing the partition's
signature so one needs to divide by $a_{k}!$ for $k=1,...,n$ which yields the
formula for the number of partitions with that signature. $\square$

The \textit{Stirling numbers} $S\left(  n,m\right)  $ \textit{of the second
kind} are the number of partitions of an $n$-element set with $m$ blocks.
Since $\sum_{k=1}^{n}a_{k}=m$ is the number of blocks, the direct formula (as
opposed to a recurrence formula) is:

\begin{center}
$S\left(  n,m\right)  =%
{\textstyle\sum\limits_{\substack{1a_{1}+2a_{2}+...+na_{n}=n\\a_{1}%
+a_{2}+...+a_{n}=m}}}
\frac{n!}{a_{1}!a_{2}!...a_{n}!\left(  1!\right)  ^{a_{1}}\left(  2!\right)
^{a_{2}}...\left(  n!\right)  ^{a_{n}}}$

Direct formula for Stirling numbers of the second kind.
\end{center}

The \textit{Bell numbers} $B\left(  n\right)  $ are the total number of
partitions on an $n$-element set so the direct formula is:

\begin{center}
$B\left(  n\right)  =\sum_{m=1}^{n}S\left(  n,m\right)  =%
{\textstyle\sum\limits_{1a_{1}+2a_{2}+...+na_{n}=n}}
\frac{n!}{a_{1}!a_{2}!...a_{n}!\left(  1!\right)  ^{a_{1}}\left(  2!\right)
^{a_{2}}...\left(  n!\right)  ^{a_{n}}}$

Direct formula for total number of partitions of an $n$-element set.
\end{center}

\subsection{The direct formulas for counting DSDs of finite vector spaces}

Each DSD $\pi=\left\{  V_{i}\right\}  _{i=1,...,m}$ of a finite vector space
of dimension $n$ also determines a number partition of $n$ using the
dimensions $n_{i}=\dim\left(  V_{i}\right)  $ in place of the set
cardinalities, and thus each DSD also has a signature $a_{1},...,a_{n}$ where
the subspaces are ordered by nondecreasing dimension and where $\sum_{k=1}%
^{n}a_{k}k=n$ and $\sum_{k=1}^{n}a_{k}=m$.

\begin{proposition}
The number of DSDs of a vector space $V$ of dimension $n$ over $GF\left(
q\right)  $ with the part-count signature $a_{1},...,a_{n}$ is:
\end{proposition}

\begin{center}
$\frac{1}{a_{1}!a_{2}!...a_{n}!}\frac{\left[  n\right]  _{q}!}{\left(  \left[
1\right]  _{q}!\right)  ^{a_{1}}...\left(  \left[  n\right]  _{q}!\right)
^{a_{n}}}q^{\frac{1}{2}\left(  n^{2}-\sum_{k}a_{k}k^{2}\right)  }$

Number of DSDs for the given signature $a_{1},...,a_{n}$ where $\sum_{k=1}%
^{n}a_{k}k=n$.
\end{center}

\noindent Proof: Reasoning first in terms of the dimensions $n_{i}$, we
calculate the number of ordered bases in a subspace of dimension $n_{1}$ of a
vector space of dimension $n$ over the finite field $GF\left(  q\right)  $
with $q$ elements. There are $q^{n}$ elements in the space so the first choice
for a basis vector is $\left(  q^{n}-1\right)  $ (excluding $0$), and since
that vector generated a subspace of dimension $q$, the choice of the second
basis vector is limited to $\left(  q^{n}-q\right)  $ elements, and so forth. Thus:

\begin{center}
$\left(  q^{n}-1\right)  \left(  q^{n}-q\right)  \left(  q^{n}-q^{2}\right)
...\left(  q^{n}-q^{n_{1}-1}\right)  $

$=\left(  q^{n}-1\right)  q^{1}\left(  q^{n-1}-1\right)  q^{2}\left(
q^{n-1}-1\right)  ...q^{n_{1}-1}\left(  q^{n-n_{1}+1}-1\right)  $

$=\left(  q^{n}-1\right)  \left(  q^{n-1}-1\right)  ...\left(  q^{n-n_{1}%
-1}-1\right)  q^{\left(  1+2+...+\left(  n_{1}-1\right)  \right)  }$

Number of ordered bases for an $n_{1}$-dimensional subspace of an
$n$-dimensional space.
\end{center}

\noindent If we then divide by the number of ordered bases for an $n_{1}%
$-dimension space:

\begin{center}
$\left(  q^{n_{1}}-1\right)  \left(  q^{n_{1}}-q\right)  ...\left(  q^{n_{1}%
}-q^{n_{1}-1}\right)  =\left(  q^{n_{1}}-1\right)  \left(  q^{n_{1}%
-1}-1\right)  ...\left(  q-1\right)  q^{\left(  1+2+...+\left(  n_{1}%
-1\right)  \right)  }$
\end{center}

\noindent we could cancel the $q^{n_{1}\left(  n_{1}-1\right)  }$ terms to
obtain the Gaussian binomial coefficient

\begin{center}
$\frac{\left(  q^{n}-1\right)  \left(  q^{n-1}-1\right)  ...\left(
q^{n-n_{1}-1}-1\right)  q^{\left(  1+2+...+\left(  n_{1}-1\right)  \right)  }%
}{\left(  q^{n_{1}}-1\right)  \left(  q^{n_{1}-1}-1\right)  ...\left(
q-1\right)  q^{\left(  1+2+...+\left(  n_{1}-1\right)  \right)  }}=\binom
{n}{n_{1}}_{q}=\frac{\left[  n\right]  _{q}!}{\left[  n-n_{1}\right]
_{q}!\left[  n_{1}\right]  _{q}!}$

Number of different $n_{1}$-dimensional subspaces of an $n$-dimensional space.
\end{center}

\noindent If instead we continue to develop the numerator by multiplying by
the number of ordered bases for an $n_{2}$-dimensional space that could be
chosen from the remaining space of dimension $n-n_{1}$ to obtain:

\begin{center}
$\left(  q^{n}-1\right)  \left(  q^{n}-q\right)  \left(  q^{n}-q^{2}\right)
...\left(  q^{n}-q^{n_{1}-1}\right)  \times\left(  q^{n}-q^{n_{1}}\right)
\left(  q^{n}-q^{n_{1}+1}\right)  ...\left(  q^{n}-q^{n_{1}+n_{2}-1}\right)  $

$=\left(  q^{n}-1\right)  \left(  q^{n-1}-1\right)  ...\left(  q^{n-n_{1}%
-n_{2}+1}-1\right)  q^{\left(  1+2+...+\left(  n_{1}+n_{2}-1\right)  \right)
}$.
\end{center}

\noindent Then dividing by the number of ordered bases of an $n_{1}%
$-dimensional space times the number of ordered bases of an $n_{2}%
$-dimensional space gives the number of different "disjoint" (i.e., only
overlap is zero subspace) subspaces of $n_{1}$-dimensional and $n_{2}%
$-dimensional subspaces.

\begin{center}
$=\frac{\left(  q^{n}-1\right)  \left(  q^{n-1}-1\right)  ...\left(
q^{n-n_{1}-n_{2}+1}-1\right)  q^{\left(  1+2+...+\left(  n_{1}+n_{2}-1\right)
\right)  }}{\left(  q^{n_{1}}-1\right)  \left(  q^{n_{1}-1}-1\right)
...\left(  q-1\right)  q^{\left(  1+2+...+\left(  n_{1}-1\right)  \right)
}\times\left(  q^{n_{2}}-1\right)  \left(  q^{n_{2}-1}-1\right)  ...\left(
q-1\right)  q^{1+2+...+\left(  n_{2}-1\right)  }}$.
\end{center}

Continuing in this fashion we arrive at the number of disjoint subspaces of
dimensions $n_{1},n_{2},...,n_{m}$ where $\sum_{i=1}^{m}n_{i}=n$:

\begin{center}
$\frac{\left(  q^{n}-1\right)  \left(  q^{n-1}-1\right)  ...\left(
q-1\right)  q^{\left(  1+2+...+\left(  n-1\right)  \right)  }}{\prod
\limits_{i=1,...,m}\left(  q^{n_{i}}-1\right)  \left(  q^{n_{i}-1}-1\right)
...\left(  q-1\right)  q^{\left(  1+2+...+\left(  n_{i}-1\right)  \right)  }%
}=\frac{\left[  n\right]  _{q}!q^{n\left(  n-1\right)  /2}}{\left[
n_{1}\right]  _{q}!q^{n_{1}\left(  n_{1}-1\right)  /2}\times...\times\left[
n_{m}\right]  _{q}!q^{n_{m}\left(  n_{m}-1\right)  /2}}$

$=\frac{\left[  n\right]  _{q}!}{\left[  n_{1}\right]  _{q}!...\left[
n_{m}\right]  _{q}!}q^{\frac{1}{2}\left[  n\left(  n-1\right)  -\sum_{i=1}%
^{m}n_{i}\left(  n_{i}-1\right)  \right]  }$.
\end{center}

\noindent There may be a number $a_{k}$ of subspaces with the same dimension,
e.g., if $n_{j}=n_{j+1}=k$, then $a_{k}=2$ so the term $\left[  n_{j}\right]
_{q}!q^{n_{j}\left(  n_{j}-1\right)  /2}\times\left[  n_{j+1}\right]
_{q}!q^{n_{j+1}\left(  n_{j+1}-1\right)  /2}$ in the denominator could be
replaced by $\left(  \left[  k\right]  _{q}!\right)  ^{a_{k}}q^{a_{k}k\left(
k-1\right)  /2}$. Hence the previous result could be rewritten in the
part-count representation:

\begin{center}
$\frac{\left[  n\right]  _{q}!}{\left(  \left[  1\right]  _{q}!\right)
^{a_{1}}...\left(  \left[  n\right]  _{q}!\right)  ^{a_{n}}}q^{\frac{1}%
{2}\left[  n\left(  n-1\right)  -\sum_{k}a_{k}k\left(  k-1\right)  \right]  }$.
\end{center}

\noindent And permuting subspaces of the same dimension $k$ yields a DSD with
the same signature, so we need to divide by $a_{k}!$ to obtain the formula:

\begin{center}
$\frac{\left[  n\right]  _{q}!}{a_{1}!...a_{n}!\left(  \left[  1\right]
_{q}!\right)  ^{a_{1}}...\left(  \left[  n\right]  _{q}!\right)  ^{a_{n}}%
}q^{\frac{1}{2}\left[  n\left(  n-1\right)  -\sum_{k}a_{k}k\left(  k-1\right)
\right]  }$.
\end{center}

\noindent The exponent on the $q$ term can be simplified since $\sum_{k}%
a_{k}k=n$:

\begin{center}
$\frac{1}{2}\left[  n\left(  n-1\right)  -\left(  \sum_{k}a_{k}k\left(
k-1\right)  \right)  \right]  =\frac{1}{2}\left[  n^{2}-n-\left(  \sum
_{k}a_{k}k^{2}-\sum_{k}a_{k}k\right)  \right]  $

$=\frac{1}{2}\left[  n^{2}-n-\left(  \sum_{k}a_{k}k^{2}-n\right)  \right]
=\frac{1}{2}\left(  n^{2}-\sum_{k}a_{k}k^{2}\right)  $.
\end{center}

This yields the final formula for the number of DSDs with the part-count
signature $a_{1},...,a_{n}$:

\begin{center}
$\frac{\left[  n\right]  _{q}!}{a_{1}!...a_{n}!\left(  \left[  1\right]
_{q}!\right)  ^{a_{1}}...\left(  \left[  n\right]  _{q}!\right)  ^{a_{n}}%
}q^{\frac{1}{2}\left(  n^{2}-\sum_{k}a_{k}k^{2}\right)  }$. $\square$
\end{center}

Note that the formula is not obtained by a simple substitution of $\left[
k\right]  _{q}!$ for $k!$ in the set partition formula due to the extra term
$q^{\frac{1}{2}\left(  n^{2}-\sum_{k}a_{k}k^{2}\right)  }$, but that it still
reduces to the classical formula for set partitions with that signature as
$q\rightarrow1$. This formula leads directly to the vector space version of
the Stirling numbers of the second kind to count the DSDs with $m$ parts and
to the vector space version of the Bell numbers to count the total number of DSDs.

Before giving those formulas, it should be noted that there is another
$q$-analog formula called "generalized Stirling numbers" (of the second
kind)--but it generalizes only one of the recurrence formulas for $S\left(
n,m\right)  $. It does not generalize the \textit{interpretation} "number of
set partitions on an $n$-element set with $m$ parts" to count the vector space
partitions (DSDs) of finite vector spaces of dimension $n$ with $m$ parts. The
Stirling numbers satisfy the recurrence formula:

\begin{center}
$S\left(  n+1,m\right)  =mS\left(  n,m\right)  +S\left(  n-1,m\right)  $ with
$S\left(  0,m\right)  =\delta_{0m}$.
\end{center}

\noindent Donald Knuth uses the braces notation for the Stirling numbers, $%
\genfrac{\{}{\}}{0pt}{}{n}{m}%
=S\left(  n,m\right)  $, and then he defines the "generalized Stirling number"
\cite[p. 436]{knuth:vol4a} $%
\genfrac{\{}{\}}{0pt}{}{n}{m}%
_{q}$ by the $q$-analog recurrence relation:

\begin{center}
$%
\genfrac{\{}{\}}{0pt}{}{n+1}{m}%
_{q}=\left(  1+q+...+q^{m-1}\right)
\genfrac{\{}{\}}{0pt}{}{n}{m}%
_{q}+%
\genfrac{\{}{\}}{0pt}{}{n}{m-1}%
_{q}$; $%
\genfrac{\{}{\}}{0pt}{}{0}{m}%
_{q}=\delta_{0m}$.
\end{center}

It is easy to generalize the direct formula for the Stirling numbers and it
generalizes the partition interpretation:

\begin{center}
$D_{q}\left(  n,m\right)  =%
{\textstyle\sum\limits_{\substack{1a_{1}+2a_{2}+...+na_{n}=n\\a_{1}%
+a_{2}+...+a_{n}=m}}}
\frac{\left[  n\right]  _{q}!}{a_{1}!...a_{n}!\left(  \left[  1\right]
_{q}!\right)  ^{a_{1}}...\left(  \left[  n\right]  _{q}!\right)  ^{a_{n}}%
}q^{\frac{1}{2}\left(  n^{2}-\sum_{k}a_{k}k^{2}\right)  }$

Number of DSDs of a finite vector space of dimension $n$ over $GF\left(
q\right)  $ with $m$ parts.
\end{center}

\noindent The number $D_{q}\left(  n,m\right)  $ is $S_{nm}$ in
\cite{stanley:expstructures}. Taking $q\rightarrow1$ yields the Stirling
numbers of the second kind, i.e., $D_{1}\left(  n,m\right)  =S\left(
n,m\right)  $. Knuth's generalized Stirling numbers $%
\genfrac{\{}{\}}{0pt}{}{n}{m}%
_{q}$ and $D_{q}\left(  n,m\right)  $ start off the same, e.g., $%
\genfrac{\{}{\}}{0pt}{}{0}{0}%
_{q}=1=D_{q}\left(  0,0\right)  $ and $%
\genfrac{\{}{\}}{0pt}{}{1}{1}%
_{q}=1=D_{q}\left(  1,1\right)  $, but then quickly diverge. For instance, all
$%
\genfrac{\{}{\}}{0pt}{}{n}{n}%
_{q}=1$ for all $n$, whereas the special case of $D_{q}(n,n)$ is the number of
DSDs of $1$-dimensional subspaces in a finite vector space of dimension $n$
over $GF\left(  q\right)  $ (see table below for $q=2$). The formula
$D_{q}\left(  n,n\right)  $ is $M\left(  n\right)  $ in \cite[Example
5.5.2(b), pp. 45-6]{stanley:ec-vol2} or \cite[Example 2.2, p. 75.]%
{stanley:expstructures}.

The number $D_{q}(n,n)$ of DSDs of $1$-dimensional subspaces is closely
related to the number of basis sets. The old formula for that number of bases
is \cite[p. 71]{lidl:finite-fields}:

\begin{center}
$\frac{1}{n!}\left(  q^{n}-1\right)  \left(  q^{n}-q\right)  ...\left(
q^{n}-q^{n-1}\right)  $

$=\frac{1}{n!}\left(  q^{n}-1\right)  \left(  q^{n-1}-1\right)  ...(q^{1}%
-1)q^{\left(  1+2+...+\left(  n-1\right)  \right)  }$

$=\frac{1}{n!}\left[  n\right]  _{q}!q^{\binom{n}{2}}\left(  q-1\right)  ^{n}$
\end{center}

\noindent\ since $\left[  k\right]  _{q}=\frac{q^{k}-1}{q-1}=1+q+q^{2}%
+...+q^{k-1}$ for $k=1,...,n$.

In the formula for $D_{q}\left(  n,n\right)  $, there is only one signature
$a_{1}=n$ and $a_{k}=0$ for $k=2,...,n$ which immediately gives the formula
for the number of DSDs with $n$ $1$-dimensional blocks and each $1$%
-dimensional block has $q-1$ choices for a basis vector so the total number of
sets of basis vectors is given by the same formula:

\begin{center}
$D_{q}\left(  n,n\right)  \left(  q-1\right)  ^{n}=\frac{\left[  n\right]
_{q}!}{a_{1}!}q^{\frac{1}{2}\left(  n^{2}-a_{1}1^{2}\right)  }\left(
q-1\right)  ^{n}=\frac{1}{n!}\left[  n\right]  _{q}!q^{\binom{n}{2}}\left(
q-1\right)  ^{n}$.
\end{center}

\noindent Note that for $q=2$, $\left(  q-1\right)  ^{n}=1$ so $D_{2}\left(
n,n\right)  $ is the number of different basis sets.

Summing the $D_{q}\left(  n,m\right)  $ for all $m$ gives the vector space
version of the Bell numbers $B\left(  n\right)  $:

\begin{center}
$D_{q}\left(  n\right)  =\sum_{m=1}^{n}D_{q}\left(  n,m\right)  =%
{\textstyle\sum\limits_{1a_{1}+2a_{2}+...+na_{n}=n}}
\frac{1}{a_{1}!a_{2}!...a_{n}!}\frac{\left[  n\right]  _{q}!}{\left(  \left[
1\right]  _{q}!\right)  ^{a_{1}}...\left(  \left[  n\right]  _{q}!\right)
^{a_{n}}}q^{\frac{1}{2}\left(  n^{2}-\sum_{k}a_{k}k^{2}\right)  }$

Number of DSDs of a vector space of dimension $n$ over $GF\left(  q\right)  $.
\end{center}

\noindent Our notation $D_{q}\left(  n\right)  $ is $D_{n}\left(  q\right)  $
in Bender and Goldman \cite{bender-goldman} and $\left\vert Q_{n}\right\vert $
in Stanley (\cite{stanley:expstructures}, \cite{stanley:ec-vol2}). Setting
$q=1$ gives the Bell numbers, i.e., $D_{1}\left(  n\right)  =B\left(
n\right)  $.

\subsection{Counting DSDs with a block containing a designated vector
$v^{\ast}$}

Set partitions have a property not shared by vector space partitions, i.e.,
DSDs. Given a designated element $u^{\ast}$ of the universe set $U$, the
element is contained in some block of every partition on $U$. But given a
nonzero vector $v^{\ast}$ in a space $V$, it is not necessarily contained in a
block of any given DSD of $V$. Some proofs of formulas use this property of
set partitions so the proofs do not generalize to DSDs.

Consider one of the formulas for the Stirling numbers of the second kind:

\begin{center}
$S\left(  n,m\right)  =\sum_{k=0}^{n-1}\binom{n-1}{k}S\left(  k,m-1\right)  $

Summation formula for $S\left(  n,m\right)  $.
\end{center}

The proof using the designated-element $u^{\ast}$ reasoning starts with the
fact that any partition of $U$ with $\left\vert U\right\vert =n$ with $m$
blocks will have one block containing $u^{\ast}$ so we then only need to count
the number of $m-1$ blocks on the subset disjoint from the block containing
$u^{\ast}$. If the block containing $u^{\ast}$ had $n-k$ elements, there are
$\binom{n-1}{k}$ blocks that could be complementary to an $\left(  n-k\right)
$-element block containing $u^{\ast}$ and each of those $k$-element blocks had
$S\left(  k,m-1\right)  $ partitions on it with $m-1$ blocks. Hence the total
number of partitions on an $n$-element set with $m$ blocks is that sum.

This reasoning can be extended to DSDs over finite vector spaces, but it only
counts the number of DSDs with a block containing a designated nonzero vector
$v^{\ast}$ (it doesn't matter which one), not all DSDs. Furthermore, it is not
a simple matter of substituting $\binom{n-1}{k}_{q}$ for $\binom{n-1}{k}$.
Each $\left(  n-k\right)  $-element subset has a unique $k$-element subset
disjoint from it (its complement), but the same does not hold in general
vector spaces. Thus given a subspace with $\left(  n-k\right)  $-dimensions,
we must compute the number of $k$-dimensional subspaces disjoint from it.

\begin{center}
Let $V$ be an $n$-dimensional vector space over $GF\left(  q\right)  $ and let
$v^{\ast}$ be a specific nonzero vector in $V$. In a DSD with an $\left(
n-k\right)  $-dimensional block containing $v^{\ast}$, how many $k$%
-dimensional subspaces are there disjoint from the $\left(  n-k\right)
$-dimensional subspace containing $v^{\ast}$? The number of ordered basis sets
for a $k$-dimensional subspace disjoint from the given $\left(  n-k\right)
$-dimensional space is:

$\left(  q^{n}-q^{n-k}\right)  \left(  q^{n}-q^{n-k+1}\right)  ...\left(
q^{n}-q^{n-1}\right)  =\left(  q^{k}-1\right)  q^{n-k}\left(  q^{k-1}%
-1\right)  q^{n-k+1}...\left(  q-1\right)  q^{n-1}$

$=\left(  q^{k}-1\right)  \left(  q^{k-1}-1\right)  ...\left(  q-1\right)
q^{\left(  n-k\right)  +\left(  n-k+1\right)  +...+\left(  n-1\right)  }$

$=\left(  q^{k}-1\right)  \left(  q^{k-1}-1\right)  ...\left(  q-1\right)
q^{k\left(  n-k\right)  +\frac{1}{2}k\left(  k-1\right)  }$
\end{center}

\noindent since we use the usual trick to evaluate twice the exponent:

\begin{center}
$\left(  n-k\right)  +\left(  n-k+1\right)  +...+\left(  n-1\right)  $

$+$\underline{$\left(  n-1\right)  +\left(  n-2\right)  +...+\left(
n-k\right)  $}

$=\left(  2n-k-1\right)  +...+\left(  2n-k-1\right)  $

$=k\left(  2n-k-1\right)  =2k\left(  k+\left(  n-k\right)  \right)
-k^{2}-k=2k\left(  n-k\right)  +k^{2}-k$.
\end{center}

Now the number of ordered basis set of a $k$-dimensional space is:

\begin{center}
$\left(  q^{k}-1\right)  \left(  q^{k-1}-1\right)  ...\left(  q-1\right)
q^{\frac{1}{2}k\left(  k-1\right)  }$
\end{center}

\noindent so dividing by that gives:

\begin{center}
$q^{k\left(  n-k\right)  }$

The number of $k$-dimensional subspaces disjoint from any $\left(  n-k\right)
$-dimensional subspace.\footnote{This was proven using M\"{o}bius inversion on
the lattice of subspaces by Crapo \cite{crapo:mobius-lattice}.}
\end{center}

\noindent Note that taking $q\rightarrow1$ yields the fact that an $\left(
n-k\right)  $-element subset of an $n$-element set has a unique $k$-element
subset disjoint from it.

Hence in the $q$-analog formula, the binomial coefficient $\binom{n-1}{k}$ is
replaced by the Gaussian binomial coefficient $\binom{n-1}{k}_{q}$ times
$q^{k\left(  n-k\right)  }$. Then the rest of the proof proceeds as usual. Let
$D_{q}^{\ast}\left(  n,m\right)  $ denote the number of DSDs of $V$ with $m$
blocks with one block containing any designated $v^{\ast}$. Then we can mimic
the proof of the formula $S\left(  n,m\right)  =\sum_{k=0}^{n-1}\binom{n-1}%
{k}S\left(  k,m-1\right)  $ to derive the following:

\begin{proposition}
Given a designated nonzero vector $v^{\ast}\in V$, the number of DSDs of $V$
with $m$ blocks one of which contains $v^{\ast}$ is:
\end{proposition}

\begin{center}
$D_{q}^{\ast}\left(  n,m\right)  =\sum_{k=0}^{n-1}$ $\binom{n-1}{k}%
_{q}q^{k\left(  n-k\right)  }D_{q}\left(  k,m-1\right)  $. $\square$
\end{center}

\noindent Note that taking $q=1$ gives the right-hand side of: $\sum
_{k=0}^{n-1}$ $\binom{n-1}{k}S\left(  k,m-1\right)  $ since $D_{1}\left(
k,m-1\right)  =S\left(  k,m-1\right)  $, and the left-hand side is the same as
$S\left(  n,m\right)  $ since \textit{every} set partition of an $n$-element
with $m$ blocks has to have a block containing some designated element
$u^{\ast}$.

Since the Bell numbers can be obtained from the Stirling numbers of the second
time as: $B\left(  n\right)  =\sum_{m=1}^{n}S\left(  n,m\right)  $, there is
clearly a similar formula for the Bell numbers:

\begin{center}
$B\left(  n\right)  =\sum_{k=0}^{n-1}\binom{n-1}{k}B\left(  k\right)  $

Summation formula for $B\left(  n\right)  $.
\end{center}

This formula can also be directly proven using the designated element
$u^{\ast}$ reasoning, so it can be similarly be extended to computing
$D_{q}^{\ast}\left(  n\right)  $, the number of DSDs of $V$ with a block
containing a designated nonzero vector $v^{\ast}$.

\begin{proposition}
Given an designated nonzero vector $v^{\ast}\in V$, the number of DSDs of $V$
with a block containing $v^{\ast}$ is:
\end{proposition}

\begin{center}
$D_{q}^{\ast}\left(  n\right)  =\sum_{k=0}^{n-1}$ $\binom{n-1}{k}%
_{q}q^{k\left(  n-k\right)  }D_{q}\left(  k\right)  $. $\square$
\end{center}

\noindent In the same manner, taking $q=1$ yields the classical summation
formula for $B\left(  n\right)  $ since $D_{1}\left(  k\right)  =B\left(
k\right)  $ and every partition has to have a block containing a designated
element $u^{\ast}$.

Furthermore the $D^{\ast}$ numbers have the expected relation:

\begin{corollary}
$D_{q}^{\ast}(n)=\sum_{m=1}^{n}D_{q}^{\ast}\left(  n,m\right)  $. $\square$
\end{corollary}

Since we also have formulas for the total number of DSDs, we have the number
of DSDs where the designated element is not in one of the blocks:
$D_{q}\left(  n,m\right)  -D_{q}^{\ast}\left(  n,m\right)  $ and $D_{q}\left(
n\right)  -D_{q}^{\ast}\left(  n\right)  $.

For a natural interpretation for the $D^{\ast}$-numbers, consider the
pedagogical model of quantum mechanics using finite-dimensional vector spaces
$V$ over $GF\left(  2\right)  =%
\mathbb{Z}
_{2}$, called "quantum mechanics over sets," QM/Sets \cite{ell:qmoversets}.
The "observables" or attributes are defined by real-valued functions on basis
sets. Given a basis set $U=\left\{  u_{1},...,u_{n}\right\}  $ for $V=%
\mathbb{Z}
_{2}^{n}\cong\wp\left(  U\right)  $, a \textit{real-valued attribute} is a
function $f:U\rightarrow%
\mathbb{R}
$. It determines a set partition $\left\{  f^{-1}\left(  r\right)  \right\}
_{r\in f\left(  U\right)  }$ on $U$ and a DSD $\left\{  \wp\left(
f^{-1}\left(  r\right)  \right)  \right\}  _{r\in f\left(  U\right)  }$ on
$\wp\left(  U\right)  $. In full QM, the important thing about an "observable"
is not the specific numerical eigenvalues, but its eigenspaces for distinct
eigenvalues, and that information is in the DSD of its eigenspaces. The
attribute $f:U\rightarrow%
\mathbb{R}
$ cannot be internalized as an operator on $\wp\left(  U\right)  \cong%
\mathbb{Z}
_{2}^{n}$ (unless its values are $0,1$), but it nevertheless determines the
DSD $\left\{  \wp\left(  f^{-1}\left(  r\right)  \right)  \right\}  _{r\in
f\left(  U\right)  }$ which is sufficient to pedagogically model many quantum
results. Hence a DSD can be thought of an "abstract attribute" (without the
eigenvalues) with its blocks serving as "eigenspaces." Then a natural question
to ask is given any nonzero vector $v^{\ast}\in V=%
\mathbb{Z}
_{2}^{n}$, how many "abstract attributes" are there where $v^{\ast}$ is an
"eigenvector"--and the answer is $D_{2}^{\ast}\left(  n\right)  $. And
$D_{2}^{\ast}\left(  n,m\right)  $ is the number of "abstract attributes" with
$m$ distinct "eigenvalues" where $v^{\ast}$ is an "eigenvector."

\section{Computing initial values for $q=2$}

In the case of $n=1,2,3$, the DSDs can be enumerated "by hand" to check the
formulas, and then the formulas can be used to compute higher values of
$D_{2}\left(  n,m\right)  $ or $D_{2}\left(  n\right)  $.

All vectors in the $n$-dimensional vector space $%
\mathbb{Z}
_{2}^{n}$ over $GF\left(  2\right)  =%
\mathbb{Z}
_{2}$ will be expressed in terms of a computational basis $\left\{  a\right\}
,\left\{  b\right\}  ,...,\left\{  c\right\}  $ so any vector $S$ in $%
\mathbb{Z}
_{2}^{n}$ would be represented in that basis as a subset $S\subseteq
U=\left\{  a,b,...,c\right\}  $. The addition of subsets (expressed in the
same basis) is the symmetric difference: for $S,T\subseteq U$,

\begin{center}
$S+T=\left(  S-T\right)  \cup\left(  T-S\right)  =\left(  S\cup T\right)
-\left(  S\cap T\right)  .$
\end{center}

\noindent Since all subspaces contain the zero element which is the empty set
$\emptyset$, it will be usually suppressed when listing the elements of a
subspace. And subsets like $\left\{  a\right\}  $ or $\left\{  a,b\right\}  $
will be denoted as just $a$ and $ab$. Thus the subspace $\left\{
\emptyset,\left\{  a\right\}  ,\left\{  b\right\}  ,\left\{  a,b\right\}
\right\}  $ is denoted as $\left\{  a,b,ab\right\}  $. A $k$-dimensional
subspace has $2^{k}$ elements so only $2^{k}-1$ are listed.

For $n=1$, there is only one nonzero subspace $\left\{  a\right\}  $, i.e.,
$\left\{  \emptyset,\left\{  a\right\}  \right\}  $, and $D_{2}\left(
1,1\right)  =D_{2}\left(  1\right)  =1$.

For $n=2$, the whole subspace is $\left\{  a,b,ab\right\}  $ and it has three
bases $\left\{  a,b\right\}  $, $\left\{  a,ab\right\}  $, and $\left\{
b,ab\right\}  $. The formula for the number of bases gives $D_{2}\left(
2,2\right)  =3$. The only $D_{2}\left(  2,1\right)  =1$ DSD is the whole space.

For $n=3$, the whole space $\left\{  a,b,c,ab,ac,bc,abc\right\}  $ is the only
$D_{2}\left(  3,1\right)  =1$ and indeed for any $n$ and $q$, $D_{q}\left(
n,1\right)  =1$. For $n=3$ and $m=3$, $D_{2}\left(  3,3\right)  $ is the
number of (unordered) bases of $%
\mathbb{Z}
_{2}^{3}$ (recall $%
\genfrac{\{}{\}}{0pt}{}{n}{n}%
_{q}=1$ for all $q$). Since we know the signature, i.e., $a_{1}=3$ and
otherwise $a_{k}=0$, we can easily compute $D_{2}\left(  3,3\right)  $:

\begin{center}
$\frac{1}{a_{1}!a_{2}!...a_{n}!}\frac{\left[  n\right]  _{q}!}{\left(  \left[
1\right]  _{q}!\right)  ^{a_{1}}...\left(  \left[  n\right]  _{q}!\right)
^{a_{n}}}q^{\frac{1}{2}\left(  n^{2}-\sum_{k}a_{k}k^{2}\right)  }$

$=\frac{1}{3!}\frac{\left[  3\right]  _{2}!}{\left(  \left[  1\right]
_{2}\right)  ^{3}}2^{\frac{1}{2}\left(  3^{2}-3\right)  }=\frac{1}{6}%
\frac{7\times3}{1}2^{\frac{1}{2}\left(  6\right)  }=28=D_{2}(3,3)$.
\end{center}

\noindent And here they are.

\begin{center}%
\begin{tabular}
[c]{|c|c|c|c|}\hline
$\left\{  a,b,c\right\}  $ & $\left\{  a,b,ac\right\}  $ & $\left\{
a,b,bc\right\}  $ & $\left\{  a,b,abc\right\}  $\\\hline
$\left\{  a,c,ab\right\}  $ & $\left\{  a,c,bc\right\}  $ & $\left\{
a,c,abc\right\}  $ & $\left\{  a,ab,ac\right\}  $\\\hline
$\left\{  a,ab,bc\right\}  $ & $\left\{  a,ab,abc\right\}  $ & $\left\{
a,ac,bc\right\}  $ & $\left\{  a,ac,abc\right\}  $\\\hline
$\left\{  b,c,ab\right\}  $ & $\left\{  b,c,ac\right\}  $ & $\left\{
b,c,abc\right\}  $ & $\left\{  b,ab,ac\right\}  $\\\hline
$\left\{  b,ab,bc\right\}  $ & $\left\{  b,ab,abc\right\}  $ & $\left\{
b,ac,bc\right\}  $ & $\left\{  b,bc,abc\right\}  $\\\hline
$\left\{  c,ab,ac\right\}  $ & $\left\{  c,ab,bc\right\}  $ & $\left\{
c,ac,bc\right\}  $ & $\left\{  c,ac,abc\right\}  $\\\hline
$\left\{  ab,ac,abc\right\}  $ & $\left\{  ab,bc,abc\right\}  $ & $\left\{
ac,bc,abc\right\}  $ & $\left\{  bc,ab,abc\right\}  $\\\hline
\end{tabular}

All bases of $%
\mathbb{Z}
_{2}^{3}$.
\end{center}

\noindent For $n=3$ and $m=2$, $D_{2}\left(  3,2\right)  $ is the number of
binary DSDs, each of which has the signature $a_{1}=a_{2}=1$ so the total
number of binary DSDs is:

\begin{center}
$D_{2}\left(  3,2\right)  =\frac{1}{1!1!}\frac{\left[  3\right]  _{2}%
!}{\left(  \left[  1\right]  !\right)  ^{1}\left(  \left[  2\right]  !\right)
^{1}}2^{\frac{1}{2}\left(  3^{2}-1-2^{2}\right)  }=\frac{7\times3}{3}%
2^{\frac{1}{2}\left(  4\right)  }=7\times4=28$.
\end{center}

\noindent And here they are:

\begin{center}%
\begin{tabular}
[c]{|c|c|c|c|}\hline
$\left\{  \left\{  a\right\}  ,\left\{  b,c,bc\right\}  \right\}  $ &
$\left\{  \left\{  a\right\}  ,\left\{  ab,ac,bc\right\}  \right\}  $ &
$\left\{  \left\{  a\right\}  ,\left\{  c,ab,abc\right\}  \right\}  $ &
$\left\{  \left\{  a\right\}  ,\left\{  b,ac,abc\right\}  \right\}  $\\\hline
$\left\{  \left\{  b\right\}  ,\left\{  a,c,ac\right\}  \right\}  $ &
$\left\{  \left\{  b\right\}  ,\left\{  ab,ac,bc\right\}  \right\}  $ &
$\left\{  \left\{  b\right\}  ,\left\{  c,ab,abc\right\}  \right\}  $ &
$\left\{  \left\{  b\right\}  ,\left\{  a,bc,abc\right\}  \right\}  $\\\hline
$\left\{  \left\{  ab\right\}  ,\left\{  b,c,bc\right\}  \right\}  $ &
$\left\{  \left\{  ab\right\}  ,\left\{  a,bc,abc\right\}  \right\}  $ &
$\left\{  \left\{  ab\right\}  ,\left\{  b,ac,abc\right\}  \right\}  $ &
$\left\{  \left\{  ab\right\}  ,\left\{  a,c,ac\right\}  \right\}  $\\\hline
$\left\{  \left\{  c\right\}  ,\left\{  a,b,ab\right\}  \right\}  $ &
$\left\{  \left\{  c\right\}  ,\left\{  ab,ac,bc\right\}  \right\}  $ &
$\left\{  \left\{  c\right\}  ,\left\{  a,bc,abc\right\}  \right\}  $ &
$\left\{  \left\{  c\right\}  ,\left\{  b,ac,abc\right\}  \right\}  $\\\hline
$\left\{  \left\{  ac\right\}  ,\left\{  a,b,ab\right\}  \right\}  $ &
$\left\{  \left\{  ac\right\}  ,\left\{  a,bc,abc\right\}  \right\}  $ &
$\left\{  \left\{  ac\right\}  ,\left\{  c,ab,abc\right\}  \right\}  $ &
$\left\{  \left\{  ac\right\}  ,\left\{  b,c,bc\right\}  \right\}  $\\\hline
$\left\{  \left\{  bc\right\}  ,\left\{  a,b,ab\right\}  \right\}  $ &
$\left\{  \left\{  bc\right\}  ,\left\{  b,ac,abc\right\}  \right\}  $ &
$\left\{  \left\{  bc\right\}  ,\left\{  c,ab,abc\right\}  \right\}  $ &
$\left\{  \left\{  bc\right\}  ,\left\{  a,c,ac\right\}  \right\}  $\\\hline
$\left\{  \left\{  abc\right\}  ,\left\{  a,b,ab\right\}  \right\}  $ &
$\left\{  \left\{  abc\right\}  ,\left\{  b,c,bc\right\}  \right\}  $ &
$\left\{  \left\{  abc\right\}  ,\left\{  a,c,ac\right\}  \right\}  $ &
$\left\{  \left\{  abc\right\}  ,\left\{  ab,ac,bc\right\}  \right\}
$\\\hline
\end{tabular}

All binary DSDs for $%
\mathbb{Z}
_{2}^{3}$.
\end{center}

\noindent The above table has been arranged to illustrate the result that any
given $k$-dimensional subspaces has $q^{k\left(  n-k\right)  }$ subspaces
disjoint from it. For $n=3$ and $k=1$, each row gives the $2^{2}=4$ subspaces
disjoint from any given $1$-dimensional subspace represented by $\left\{
a\right\}  $, $\left\{  b\right\}  $,..., $\left\{  abc\right\}  $. For
instance, the four subspaces disjoint from the subspace $\left\{  ab\right\}
$ (shorthand for $\left\{  \emptyset,\left\{  a,b\right\}  \right\}  $) are
given in the third row since those are the "complementary" subspaces that
together with $\left\{  ab\right\}  $ form a DSD.

For $q=2$, the initial values up to $n=6$ of $D_{2}\left(  n,m\right)  $are
given the following table.

\begin{center}%
\begin{tabular}
[c]{|c||c|c|c|c|c|c|c|}\hline
$n\backslash m$ & $0$ & $1$ & $2$ & $3$ & $4$ & $5$ & $6$\\\hline\hline
$0$ & $1$ &  &  &  &  &  & \\\hline
$1$ & $0$ & $1$ &  &  &  &  & \\\hline
$2$ & $0$ & $1$ & $3$ &  &  &  & \\\hline
$3$ & $0$ & $1$ & $28$ & $28$ &  &  & \\\hline
$4$ & $0$ & $1$ & $400$ & $1,680$ & $840$ &  & \\\hline
$5$ & $0$ & $1$ & $10,416$ & $168,640$ & $277,760$ & $83,328$ & \\\hline
$6$ & $0$ & $1$ & $525,792$ & $36,053,248$ & $159,989,760$ & $139,991,040$ &
$27,998,208$\\\hline
\end{tabular}

$D_{2}\left(  n,m\right)  $ with $n,m=1,2,...,6$.
\end{center}

The seventh row $D_{2}\left(  7,m\right)  $ for $m=0,1,...,7$ is: $0$, $1$,
$51116992$, $17811244032$, $209056841728$, $419919790080$, $227569434624$, and
$32509919232$ which sum to $D_{2}\left(  7\right)  $.

The row sums give the values of $D_{2}\left(  n\right)  $ for $n=0,1,2,...,7$.

\begin{center}%
\begin{tabular}
[c]{|c||c|}\hline
$n$ & $D_{2}\left(  n\right)  $\\\hline\hline
$0$ & $1$\\\hline
$1$ & $1$\\\hline
$2$ & $4$\\\hline
$3$ & $57$\\\hline
$4$ & $2,921$\\\hline
$5$ & $540,145$\\\hline
$6$ & $364,558,049$\\\hline
$7$ & $906,918,346,689$\\\hline
\end{tabular}

$D_{2}\left(  n\right)  $ for $n=0,1,...,7$.
\end{center}

We can also compute the $D^{\ast}$ examples of DSDs with a block containing a
designated element. For $q=2$, the $D_{2}^{\ast}\left(  n,m\right)  $ numbers
for $n,m=0,1,...,7$ are given in the following table.

\begin{center}

\begin{tabular}
[c]{|c||c|c|c|c|c|c|c|c|}\hline
$n\backslash m$ & $0$ & $1$ & $2$ & $3$ & $4$ & $5$ & $6$ & $7$\\\hline\hline
$0$ & $1$ &  &  &  &  &  &  & \\\hline
$1$ & $0$ & $1$ &  &  &  &  &  & \\\hline
$2$ & $0$ & $1$ & $2$ &  &  &  &  & \\\hline
$3$ & $0$ & $1$ & $16$ & $12$ &  &  &  & \\\hline
$4$ & $0$ & $1$ & $176$ & $560$ & $224$ &  &  & \\\hline
$5$ & $0$ & $1$ & $3456$ & $40000$ & $53760$ & $13440$ &  & \\\hline
$6$ & $0$ & $1$ & $128000$ & $5848832$ & $20951040$ & $15554560$ & $2666496$ &
\\\hline
$7$ & $0$ & $1$ & $9115648$ & $1934195712$ & $17826414592$ & $30398054400$ &
$14335082496$ & $1791885312$\\\hline
\end{tabular}

Number of DSDs $D_{2}^{\ast}\left(  n,m\right)  $ containing any given nonzero
vector $v^{\ast}$
\end{center}

For $n=3$ and $m=2$, the table says there are $D_{2}^{\ast}\left(  3,2\right)
=16$ DSDs with $2$ blocks one of which contains a given vector, say $v^{\ast
}=ab$ which represents $\left\{  a,b\right\}  $, and here they are.

\begin{center}%
\begin{tabular}
[c]{|c|c|c|c|}\hline
$\left\{  \left\{  ab\right\}  ,\left\{  b,c,bc\right\}  \right\}  $ &
$\left\{  \left\{  ab\right\}  ,\left\{  a,bc,abc\right\}  \right\}  $ &
$\left\{  \left\{  ab\right\}  ,\left\{  b,ac,abc\right\}  \right\}  $ &
$\left\{  \left\{  ab\right\}  ,\left\{  a,c,ac\right\}  \right\}  $\\\hline
$\left\{  \left\{  c\right\}  ,\left\{  a,b,ab\right\}  \right\}  $ &
$\left\{  \left\{  c\right\}  ,\left\{  ab,ac,bc\right\}  \right\}  $ &
$\left\{  \left\{  ac\right\}  ,\left\{  c,ab,abc\right\}  \right\}  $ &
$\left\{  \left\{  bc\right\}  ,\left\{  c,ab,abc\right\}  \right\}  $\\\hline
$\left\{  \left\{  ac\right\}  ,\left\{  a,b,ab\right\}  \right\}  $ &
$\left\{  \left\{  a\right\}  ,\left\{  ab,ac,bc\right\}  \right\}  $ &
$\left\{  \left\{  a\right\}  ,\left\{  c,ab,abc\right\}  \right\}  $ &
$\left\{  \left\{  abc\right\}  ,\left\{  a,b,ab\right\}  \right\}  $\\\hline
$\left\{  \left\{  bc\right\}  ,\left\{  a,b,ab\right\}  \right\}  $ &
$\left\{  \left\{  b\right\}  ,\left\{  ab,ac,bc\right\}  \right\}  $ &
$\left\{  \left\{  b\right\}  ,\left\{  c,ab,abc\right\}  \right\}  $ &
$\left\{  \left\{  abc\right\}  ,\left\{  ab,ac,bc\right\}  \right\}
$\\\hline
\end{tabular}

Two-block DSDs of $\wp\left(  \left\{  a,b,c\right\}  \right)  $ with a block
containing $ab=\left\{  a,b\right\}  $.
\end{center}

\noindent The table also says there are $D_{2}^{\ast}\left(  3,3\right)  =12$
basis sets containing \textit{any} given element which we could take to be
$v^{\ast}=abc=\left\{  a,b,c\right\}  $, and here they are.

\begin{center}%
\begin{tabular}
[c]{|c|c|c|c|}\hline
$\left\{  a,b,abc\right\}  $ & $\left\{  b,ab,abc\right\}  $ & $\left\{
a,c,abc\right\}  $ & $\left\{  b,bc,abc\right\}  $\\\hline
$\left\{  a,ab,abc\right\}  $ & $\left\{  a,ac,abc\right\}  $ & $\left\{
b,c,abc\right\}  $ & $\left\{  c,ac,abc\right\}  $\\\hline
$\left\{  ab,ac,abc\right\}  $ & $\left\{  ab,bc,abc\right\}  $ & $\left\{
ac,bc,abc\right\}  $ & $\left\{  bc,ab,abc\right\}  $\\\hline
\end{tabular}

Three-block DSDs (basis sets) of $\wp\left(  \left\{  a,b,c\right\}  \right)
$ with a basis element $abc=\left\{  a,b,c\right\}  $.
\end{center}

Summing the rows in the $D_{2}^{\ast}\left(  n,m\right)  $ table gives the
values for $D_{2}^{\ast}\left(  n\right)  $ for $n=0,1,...,7$.

\begin{center}%
\begin{tabular}
[c]{|c||c|}\hline
$n$ & $D_{2}^{\ast}\left(  n\right)  $\\\hline\hline
$0$ & $1$\\\hline
$1$ & $1$\\\hline
$2$ & $3$\\\hline
$3$ & $29$\\\hline
$4$ & $961$\\\hline
$5$ & $110,657$\\\hline
$6$ & $45,148,929$\\\hline
$7$ & $66,294,748,161$\\\hline
\end{tabular}

$D_{2}^{\ast}\left(  n\right)  $ for $n=0,1,...,7$.
\end{center}

The integer sequence $D_{2}\left(  n,n\right)  $ for $n=0,1,2,...$ is known
as: $A053601$"Number of bases of an $n$-dimensional vector space over $GF(2)$"
in the \textit{On-Line Encyclopedia of Integer Sequences} (https://oeis.org/).
In the $D_{2}\left(  n,m\right)  $ table, that sequence is on the top-right
edge of the triangle: $1,1,3,28,840,...$. But the (flattened) general sequence
for $D_{2}\left(  n,m\right)  $, the number of direct-sum decompositions with
$m$ parts (generalized Stirling numbers of the second kind for finite vector
spaces) and the sequence $D_{2}\left(  n\right)  $ for the total number of
direct-sum decompositions of an $n$-dimensional vector space over $GF(2)$
(generalized Bell numbers for finite vector spaces), as well as their
$D^{\ast}$ counterparts, do not now (March 2016) appear in the
\textit{Encyclopedia}.

\end{document}